\documentclass[12pt]{amsart}
\usepackage{amsfonts}
\usepackage{amsmath,amscd,lscape}
\usepackage{amsthm}
\usepackage{amssymb}
\usepackage{latexsym}
\usepackage{mathrsfs}
\usepackage{tikz}

\usepackage{color,ulem}

\newtheorem{prop}{Proposition}[section]

\newtheorem{defn}{Definition}[section]

\newtheorem{rem}{Remark}

\newcommand{\bea}{\begin{eqnarray}}
\newcommand{\eea}{\end{eqnarray}}
\newcommand{\beano}{\begin{eqnarray*}}
\newcommand{\eeano}{\end{eqnarray*}}
\newcommand{\beq}{\begin{equation}}
\newcommand{\eeq}{\end{equation}}

\numberwithin{equation}{section}

\def\su{{\mathfrak{su}(2)}}
\def\sun{{\mathfrak{su}(1,1)}}

\def\cA{{\mathcal A}}        
        
    \def\cH{{\mathcal H}}

\numberwithin{equation}{section}
\begin{document}

%\title[On Heun algebras of the Krawtchouk and Meixner types]{On Heun algebras of the Krawtchouk and Meixner types}
\title[ Heun algebras of Lie type]{Heun algebras of Lie type}
%\dedicatory{}

\author[N. Cramp\'e]{Nicolas Cramp\'e$^{\dagger,*}$}
\address{$^\dagger$ Institut Denis-Poisson CNRS/UMR 7013 - Universit\'e de Tours - Universit\'e d'Orl\'eans; 
Parc de Grammont, 37200 Tours, France.}
\email{crampe1977@gmail.com}

\author[L. Vinet]{Luc Vinet$^{*}$}
\address{$^*$ Centre de recherches, math\'ematiques, Universit\'e de Montr\'eal; 
P.O. Box 6128, Centre-ville Station;
Montr\'eal (Qu\'ebec), H3C 3J7, Canada.}
\email{vinet@CRM.UMontreal.ca}

\author[A. Zhedanov]{Alexei Zhedanov$^{\ddagger}$}
\address{$^{\ddagger}$ Department of Mathematics, School of Information;
Renmin University of China;
Beijing 100872, China.}
\email{zhedanov@yahoo.com}

\begin{abstract}
We introduce Heun algebras of Lie type. They are obtained from bispectral pairs associated to simple or solvable Lie algebras of  dimension three or four.
For $\mathfrak{su}(2)$, this leads to the Heun-Krawtchouk algebra. The corresponding Heun-Krawtchouk operator is identified as the Hamiltonian of the quantum 
analogue of the Zhukovski-Voltera gyrostat. For $\mathfrak{su}(1,1)$, one obtains the Heun algebras attached to the Meixner, Meixner-Pollaczek and Laguerre polynomials.
These Heun algebras are shown to be isomorphic the the Hahn algebra. Focusing on the harmonic oscillator algebra $\mathfrak{ho}$ leads to the Heun-Charlier algebra. The connections to orthogonal polynomials are achieved through realizations of the underlying Lie algebras in terms of difference and differential operators. In the $\mathfrak{su}(1,1)$ cases, it is observed that the Heun operator can be transformed into the Hahn, Continuous Hahn and Confluent Heun operators respectively. 
\end{abstract}

\maketitle

%{\small MSC:\ 81R50;\ 81R10;\ 81U15.}
%

%{{\small  {\it \bf Keywords}: }}
%
\vspace{0cm}

\vspace{3mm} 
 
\section{Introduction}
This paper elaborates on Heun algebras. This subject is associated to the notion of algebraic Heun operator recently introduced \cite{GVZ_band}. Let us recall the main points. It is known that the properties of the orthogonal polynomials of the Askey scheme can be encoded algebraically. Indeed the recurrence operator X and Y, the one of which the polynomials are eigenfunctions, form a bispectral pair that generate a quadratic algebra. For the Askey-Wilson (AW) polynomials sitting at the top of the q-tableau, the defining relations of the AW algebra read \cite{Z91}:
\begin{eqnarray} 
{}[X,[X,Y]]= \rho XYX + a_1 X^2 + a_2 \{X,Y\} + a_3 X + a_4 Y + a_5\label{eq:AW1} \\
{}[Y,[Y,X]]= \rho YXY + a_1\{X,Y\} + a_2 Y^2 + a_3 Y + a_6 X + a_7\label{eq:AW2}
\end{eqnarray}
where $\rho=q^2+q^{-2}-2$ and ${a_i},  i=1,...7$ are real parameters. When $q =1$ and hence  $\rho = 0$, the cubic terms in \eqref{eq:AW1} and \eqref{eq:AW2} drop and the reduced relations define the Racah algebra \cite{GZ}, \cite{GVZ}. To all families of hypergeometric polynomials belonging to the Askey classification, there corresponds an algebra of that type which is a special case or a limit of the AW algebra.\\
The algebraic Heun operator is the generic bilinear element constructed from the generators $X,Y$
\begin{equation}
 W= r_1[X , Y]+r_2 \{X , Y\}+r_3X+r_4Y+r_5 \label{eq:Wintro}
\end{equation}
%\begin{equation}
%W= \tau_1 XY + \tau_2 YX + \tau_3 X + \tau_4 Y + \tau_0 
%\end{equation}
where $r_i, \: i=1,2,\dots,5$ are arbitrary constants. By Heun algebras, we mean the ones generated by the pairs $(X, W)$ or $(Y, W)$. 
The kernel of $W$ can be viewed as solutions of an ordinary eigenvalue problem if the parameters $r_i, i=1, \dots, 4$ are fixed and $r_5$ is the eigenvalue, or of a generalized problem if $W$ is regarded as a multiparameter linear pencil.
When X is multiplication by the variable and Y the hypergeometric operator, the algebra that is realized is the Jacobi one where $ \rho = a_1 = a_6 = a_7 = 0$. It has been shown \cite{GVZ2} that in this case $W$ coincides with the differential Heun operator and the equation $W\psi = 0$ amounts to the standard Heun equation with four regular Fuchsian singularities. The name algebraic Heun operator has been coined as a result of this observation; as already explained this concept applies to all bispectral problems and associates an operator of Heun type to the polynomials of the Askey scheme. One would then refer to the Heun operator say of Racah type and correspondingly to the Heun-Racah algebra for example.
The exploration of these structures has been initiated recently with a focus on the   \textquotedblleft higher" polynomials. Attention was first paid to the Hahn polynomials \cite{VZ_HH}; this led naturally to a finite difference version of the Heun equation on the uniform lattice. Examining the little and big q-Jacobi polynomials from this angle allowed \cite{BVZ_Heun} to give context to q-Heun operators that had been identified \cite{Takemura} in connection with integrable models. Last, the Heun-Askey-Wilson algebra was thoroughly examined \cite{BTVZ}.
We here wish to look in a similar way at the \textquotedblleft lower" polynomials and study the Heun structures attached to bispectral pairs that generate a Lie algebra. This will lead to a description of the Heun algebras associated to the Krawtchouk, Meixner, Meixner-Pollaczek, Laguerre and Charlier polynomials while establishing and clarifying the general foundations of the subject.\\
The presentation will unfold as follows. We shall begin in Section 2 with a definition of the Heun algebra of Lie type. It will be shown to be generically isomorphic to the Hahn algebra \cite{GLZ,GVZ, FGVVZ} which is obtained from the AW algebra by setting $\rho = 0$ as well as $a_2 = 0$ in (\ref{eq:AW1}) and (\ref{eq:AW2}).  We shall pursue by introducing the Lie algebra for a bispectral pair X and Y with the so-called \textquotedblleft ladder" property, this will amount to dropping the nonlinear terms in the AW relations. The corresponding W will then be seen to generate together with $X$ (or $Y$), the Heun algebra of Lie type. The definitions of the Lie algebras $\mathfrak{su}(2)$, $\mathfrak{su}(1,1)$ and $\mathfrak{ho}$ - the harmonic oscillator one, will be recorded at the end of the section. Section 3 will focus on the case where the Lie algebra for the bispectral pair is $\mathfrak{su}(2)$. It will specify the associated Heun algebra and will indicate that it cannot be mapped to the Hahn algebra over $\mathbb{R}$. It will explain in addition that the Heun operator can be viewed as as the Hamiltonian of the quantum 
analog of the Zhukovski-Voltera gyrostat. The case of $\mathfrak{su}(1,1)$ will be treated in Section 4. Three situations will be distinguished depending on whether the generator $X$ is of elliptic, hyperbolic or parabolic type. The Heun algebras will again be characterized and seen in these instances to be isomorphic to the Hahn algebra. Section 5 will be dedicated to the harmonic oscillator algebra $\mathfrak{ho}$; the associated Heun algebra will be determined and observed not to be equivalent to the Hahn algebra. By recalling in Section 6 certain realizations  of $\mathfrak{su}(2)$, $
\mathfrak{su}(1,1)$ and $\mathfrak{ho}$ in terms of difference and differential operators, we shall recognize in each case that $X$ and $Y$  become the bispectral operators of the Krawtchouk, Meixner, Meixner-Pollaczek, Laguerre and Charlier polynomials confirming that the Heun algebras of the Lie type identified  are to be associated to each of these families of orthogonal polynomials. In view of the fact that the Heun-Meixner, Heun-Meixner-Pollaczek and Heun-Laguerre algebras are isomorphic to the Hahn algebra,  it is further shown that the corresponding Heun operators can be transformed into the Hahn, Continuous Hahn and Confluent Heun operators respectively. The paper will end with concluding remarks in Section 7.

%\section{Heun algebra of the $sl_2$ type}
\section{Heun algebra of the Lie type}
\begin{defn} The Heun algebra $\cH$ of the Lie type is generated by $X$ and $W$ with the following defining relations 
\begin{eqnarray}
{} [[X,W],X]&=& x_0 +x_1 X+x_2X^2+x_3W\ ,\label{eq:H1}\\
{} [W,[X,W]]&=& y_0 + y_1X +y_2X^2+y_3X^3+x_1W+x_2\{X,W\} \label{eq:H2} ,
\end{eqnarray}
where $x_i$ and $y_i$ ($i=0,1,2,3$) are free parameters.
\end{defn}
In the Heun algebra $\cH$, the following element
\begin{equation}
 \Omega=z_1X+z_2W+z_3\{X,W\}+z_4XWX +z_5X^2+z_6W^2+z_7([X,W])^2+z_{8}X^3+z_{9}X^4
\end{equation}
is central if the parameters $z_i$ are given by
\begin{align}
 &z_1= 2y_0 - x_3y_2/3\ , && z_2=- x_2x_3 + 2x_0\ , &&  z_3=x_1\ ,\\
 &  z_4=2x_2\ , && z_5= y_1 -x_3y_3/2\ , &&z_6=x_3\ , \\
 & z_7=1\ , &&  z_{8}= 2y_2/3\ , &&z_{9}=y_3/2\ .
\end{align}

\begin{prop} \label{propH}Generically, the Heun algebra of Lie type is isomorphic to the Hahn algebra.
\end{prop}
 To show that, we introduce the following invertible map between the pairs ($X$, $W$) and ($X$, $\overline W$) where
 \begin{equation}
  \overline W = W +\mu X+\nu X^2
 \end{equation}
with $\mu$ and $\nu$ solutions of
\begin{equation}
 2\nu^2x_3-4\nu x_2+y_3=0\quad\text{and} \qquad y_2-3\mu x_2 -3\nu x_1 +3\mu\nu x_3=0\ . \label{eq:mn}
\end{equation}
It follows that $X$ and $\overline W$ satisfy the Hahn algebra 
\begin{eqnarray}
 {} [[X,\overline W],X]&=& \overline x_0 +\overline x_1 X+\overline x_2X^2+\overline x_3\overline W\\
{} [\overline W,[X,\overline W]]&=& \overline y_0 + \overline y_1X +\overline x_1\overline W+\overline x_2\{X,\overline W\}
\end{eqnarray}
where
\begin{eqnarray}
 && \overline x_0 =x_0 \ , \qquad \overline x_1=x_1-\mu x_3 \ , \qquad \overline x_2=x_2-\nu x_3 \ , \qquad \overline x_3=x_3 \label{eq:c1}\\
 &&\overline y_0=  y_0 -\mu x_0 \ , \qquad \overline y_1=y_1-2\mu x_1-2\nu x_0+\mu^2x_3 \label{eq:c2}
\end{eqnarray}
\begin{rem}
The word generically is used in the statement of Proposition \ref{propH} to indicate that it assumes freeness of the parameters. We shall observe
in the specific cases that will be discussed in the following that there are instances where the parameters need to belong to $\mathbb C$
or where there are no solutions if the parameters are not constrained. In such situations, the algebraic equivalence would not prevail.
\end{rem}
\begin{rem}

It is interesting to point out that the truncated reflection algebra attached to the Yangian of sl(2) has been shown in  \cite{CRVZ} to be isomorphic to the Hahn algebra defined above.
\end{rem}
%\subsection{Heun operators for rank one Lie algebras}

We bring at this point the generic Lie algebra $\cA$ for a bispectral pair $X$ and $Y$. Start with the $\it{linear}$  relations:
\begin{eqnarray}
&& [X , Y] = Z\, \\
&& [Z , X] = aX + c_2Y + d_2\, \label{eq:ZX} \\
&& [Y, Z] = bY + c_1X + d_1\label{eq:YZ}, 
\end{eqnarray}
where $a$, $b$, $c_1$, $c_2$, $d_1$ and $d_2$ are real parameters. With the presence of the constants $d_2$ and $d_1$ in the right-hand sides of \eqref{eq:ZX} and \eqref{eq:YZ}, we are de facto assuming in general the presence of an additional central element $I$ which we will omit writing. It is readily seen that one must have $a = b$ for the Jacobi identity $ [X , [Y , Z]] + [Y , [Z ,  X]] + [Z , [X , Y]] = 0$ to be verified. Now assume that $X$ has a nondegenerate discrete real spectrum $ \{\lambda_p,  \ p \in \mathbb Z\} $ in a representation space with $u_p$ the corresponding eigenvectors. It is easy to see that for $Y$ to act tridiagonally on that space: $Yu_p=f_{p+1}u_{p+1}+g_pu_p+f_pu_{p-1}$, in view of \eqref{eq:ZX}, 
one must have
\begin{equation}
-(\lambda_{p+1}-\lambda_p)^2=c_2.
\end{equation}
This implies the condition $c_2 < 0$. To have (Leonard) duality requires that $X$ and $Y$ be on the same footing. To that end, one would thus ask that there be reciprocally a basis
in which $Y$ has also a real discrete spectrum and $X$ is tridiagonal. This would necessitate in turn that $c_1 < 0$. In the following, we shall want however to cover more
general bispectral situations where the spectrum of $Y$ could be pure imaginary or continuous. We shall hence not impose restrictions on $c_1$. As a consequence, the Heun algebras to be defined will relate to bispectral problems that are not only the purely discrete ones.
\begin{defn} 
\label{defn1}
The generic Lie algebra $\cA$ over $\mathbb{R}$ for the bispectral pair $X$ and $Y$ is defined by \eqref{eq:ZX} and \eqref{eq:YZ} with 
\begin{equation}
a = b,  \qquad  c_2 < 0.
\end{equation}
As already indicated this corresponds to eliminating the non-linear terms in \eqref{eq:AW1} and \eqref{eq:AW2}. \
\end{defn}

The following element
\begin{equation}
 C= 2d_1 X+2d_2 Y+b \{X, Y\} +c_1 X^2+c_2 Y^2+ Z^2
\end{equation}
is central in $\cA$.

\begin{rem}\label{rem1}
It will be useful as we proceed to take note of the following two obvious facts. First, like the equations \eqref{eq:AW1} and \eqref{eq:AW2} of which they are a special case,
the equations  \eqref{eq:ZX} and \eqref{eq:YZ} preserve their form under affine transformations of $X$ and $Y$. Second, it is also clear that the elements $\underline X=VXV^{-1}$
and $\underline Y=VYV^{-1}$ obtained from performing on $X$ and $Y$ the same similarity operation will satisfy the same relations as $X$ and $Y$.
\end{rem}

\begin{prop}
\label{prop1}
Let $X$ and $Y$ be the generators of $\cA$. The Heun operator
\begin{equation}
 W= r_1[X , Y]+r_2 \{X , Y\}+r_3X+r_4Y+r_5 \label{eq:W}
\end{equation}
together with $X$ satisfy the relations of the Heun algebra $\cH$ of the Lie type with
\begin{eqnarray}
 &&\; \:  \:  \: x_0=d_2r_4-c_2r_5   \qquad x_1=2d_2r_2 - c_2r_3 + br_4 \qquad x_2=2br_2 \qquad x_3=c_2\\
 && \:  \:  \:  \;y_0=- d_2r_4r_3 - 2d_2r_5r_2 + c_2r_5r_3 + d_1r_4^2 + (d_1c_2-bd_2)(r_1^2 -r_2^2) - br_5r_4-2r_2r_4 C\\
 &&\:  \:  \:   \;y_1=- 4d_2r_3r_2 + c_2r_3^2 + 8d_1r_4r_2 + c_1r_4^2 + (c_1c_2- b^2)(r_1^2 -r_2^2) - 2br_4r_3 - 4br_5r_2 -4r_2^2C\\
 && \:  \:  \:  \; y_2=12d_1r_2^2 + 6c_1r_4r_2 - 6br_3r_2\qquad y_3=8c_1r_2^2.
\end{eqnarray}
\end{prop}
\begin{rem}
The defining relations of $\cA$ are invariant under the exchange of $X$ and $Y$ provided one performs: $c_1 \leftrightarrow c_2$, $d_1 \leftrightarrow d_2$.
Taking into account the definition of $W$ \eqref{eq:W}, we see that the algebra realized by $Y$ and $W$ is the same as the one stemming from the pair $X$, $Y$ with the coefficients  obtained from the ones given in Proposition \ref{prop1}  by doing also the changes $r_1 \rightarrow -r_1$, $r_3 \leftrightarrow r_4$.
\end{rem}
\begin{rem}
The central element $\Omega$ of $\cH$ with the parameters given above is 
\begin{equation}
 \Omega=(  r_4^2+c_2(r_1^2 - r_2^2))C- 2c_2r_5r_2b - 2c_2r_5^2 + 2r_4c_2r_2d_1 + 2d_2r_4r_5 + d_2^2(r_1^2 - r_2^2)\;.
\end{equation}
\end{rem}

In the next three sections, we shall identify the specific Heun algebras that are obtained when the algebra $\cA$ is isomorphic to $\su$, $\sun$ 
or the oscillator algebra $\mathfrak{ho}$. Let us record here the definitions of these algebras.\\

The Lie algebras $\su$ and $\sun$ are generated by $J_1$, $J_2$ and $J_3$ subject to
\begin{equation}
[J_1,J_2]=\pm iJ_3,\  \qquad [J_2,J_3]=iJ_1, \ \qquad [J_3,J_1]=iJ_2.\label{eq:crsu}
\end{equation}
The upper (resp. lower) sign is associated to $\su$  (resp. $\sun$). The quadratic Casimir of these algebras is
\begin{equation}
 \mathfrak{c}=   J_1^2 +J_2^2 \pm J_3^2 \;.\label{eq:csu}
\end{equation}

The harmonic oscillator algebra $\mathfrak{ho}$  is solvable and has four generators $N$, $A$, $A^{\dag}$ and $I$ that satisfy the commutation relations \cite{Streater}:
\begin{equation}
[N,A]=-A, \qquad [N,A^{\dag}]=A^{\dag},  \qquad [A,A^{\dag}]=I, \qquad [I,N]=[I,A]=[I,A^{\dag}]=0.\label{eq:reho}
\end{equation}
It is familiar that $\mathfrak{k}=N-A^{\dag}A$ is central in $\mathfrak{ho}$.

\section{The case of the Lie algebra $\su$}

We shall focus here on $\su$ and assume therefore in this section that the generators $(J_1, J_2, J_3)$ obey the commutation relations
of this algebra, that is those corresponding to the upper signs in \eqref{eq:crsu} and \eqref{eq:csu}.

\subsection{Specialization\label{sec:spec2}}

Take $X$ and $Y$ to be 
\begin{equation}
 X=\alpha J_3 +\beta J_1 \quad, \quad Y= J_3 \label{eq:s}
\end{equation}
with the coefficients $\alpha=\cos\theta$ and $\beta=-\sin\theta$. 

Let us first explain that there is no loss of generality with this choice keeping in mind Remark \ref{rem1}.
We can always pick $Y$ in a preferred direction and allowing for scaling have it as above. 
Assume then that $X$ is initially of the generic form $ l_1J_1+l_2J_2+l_3J_3$. The adjoint representation will 
transform this $\su$ element while preserving $l^2=l_1^2+l_2^2+l_3^2$. Under a rotation about the 3-axis, $Y$ does not change
and the other one can be transformed into $\alpha J_3 +\beta J_1$ with $\alpha^2+\beta^2 =l^2$. Scaling finally by $1/l$ we arrive at the expression for $X$ in \eqref{eq:s}.

It is readily checked that $(X,Y)$ satisfy the relations \eqref{eq:ZX} and \eqref{eq:YZ} of the algebra $\cA$ with
\begin{equation}
 d_1=d_2=0\quad , \qquad c_1=-1\quad, \quad b=\alpha \quad\text{and}\qquad c_2=-(\alpha^2+\beta^2)=-1\ .\label{eq:par}
\end{equation}
The fact that $d_1=d_2=0$ indicates that there is no center and that the algebra is three-dimensional.
The Casimir $C$ of $\cA$ is related to $\mathfrak{c}$  as follows:
\begin{equation}
 C=\alpha \{X,Y\}-X^2+c_2 Y^2+([X,Y])^2 = - \beta^2 \mathfrak{c}.\label{eq:Cas}
\end{equation}

Observe that $X$ is obtained from $Y$ through an automorphism:
\begin{equation}
X=U(\theta) Y U(\theta)^{-1} \quad \text{with}\qquad  U(\theta)= e^{i\theta J_2}. \label{eq:conj}
\end{equation}
Let $W$ be given by \eqref{eq:W}. From the results of the previous section, we can see that the algebras generated by the pairs $(X, W)$ and $(Y,W)$ will be similar
We shall therefore only focus on the Heun algebra generated by $X$ and $W$.

\begin{prop}\label{propHeun} The elements $X$ and $W$ associated to $\su$ as per \eqref{eq:s} satisfy the relations of the Heun algebra of Lie type with
\begin{eqnarray}
&&x_0=-c_2 r_5\ ,   \qquad x_1=r_4\alpha-c_2 r_3\ , \qquad x_2=2r_2\alpha\ , \qquad x_3=c_2\ ,\\ \label{eq:Ha1}
&&y_0=r_5(c_2r_3-\alpha r_4)-2r_2r_4 C\ ,\\
&&y_1=(r_2^2-r_1^2)(c_2+\alpha^2)-4r_2r_5\alpha+c_2 r_3^2-2\alpha r_3r_4-r_4^2 -4r_2^2C\ ,\\
&&y_2=-6r_2(r_4+\alpha r_3)\ , \qquad y_3=-8 r_2^2\ ,\label{eq:Ha2}
\end{eqnarray}
with $\alpha=\cos\theta$ and $c_2 =-(\alpha^2+\beta^2)=-1$.
\end{prop}
Let us remark that the central element $\Omega$ of $\cH$ with the parameters given above is 
\begin{equation}
 \Omega=(  c_2(r_1^2-r_2^2)+r_4^2 )C-c_2r_5(2r_2\alpha +r_5)\label{eq:omega}  \;.
\end{equation}

\begin{prop}
The Heun algebra of $\su$ type is not isomorphic to the Hahn algebra over $\mathbb{R}$.
\end{prop}

This is seen by observing that relation \eqref{eq:mn} in Proposition \ref{propH} for the parameters $\mu$ and $\nu$ has necessarily complex solutions in this $\su$ case, namely
\begin{equation}
 \nu=-2r_2 \exp(\pm i\theta)\quad, \qquad \mu=-r_3-\exp(\pm i\theta)r_4\ .
\end{equation}
This is indicative of the fact that the dual Hahn polynomials cannot be obtained by tridiagonalization of the Krawtchouk 
difference operator (see later).

\subsection{$\su$ Heun operators and Hamiltonians for generalized tops}

We have already noted (see Remark \ref{rem1}) that the bispectral pair of operators $X,Y$ can be chosen in various equivalent ways from the algebraic perspective. 
Different choices will however lead to modified expressions for the corresponding algebraic Heun operator. To point out the connection between Heun operators
of $ \su$ type and quantum Hamiltonians for tops, instead of \eqref{eq:s}, it will be convenient to rather adopt for $X$ and $Y$
\begin{equation}
X= J_3 + \beta J_1, \qquad Y = J_3 -\beta J_1\ , \label{eq:XY_sym} 
\end{equation}
with $\beta \ne 0$ an arbitrary parameter. One easily convinces oneself that this choice is equivalent to the preceding one. Start from $X$ and $Y$ as in  \eqref{eq:s} and conjugate
these operators by $U(-\theta/2)$ with $U(\theta)$ given in \eqref{eq:conj}. Scaling by $\sec(\theta/2)$ then gives the $X$ and $Y$ of \eqref{eq:XY_sym}.
It is easy to see that these two operators also satisfy the relations   \eqref{eq:ZX}-\eqref{eq:YZ} with
\begin{equation}
 d_1=d_2=0,\quad   \quad a=b=(1-\beta^2) \qquad\text{and}\qquad c_1=c_2=-(1+\beta^2). \ 
\end{equation}
One finds that the algebraic Heun operator \eqref{eq:W} is in this case of the form
\begin{equation}
W = \sigma \left( J_3^2 - \beta^2 J_1^2 \right) +m_1 J_1 + m_2 J_2 +  m_3 J_3+ m_4 \, \label{eq:W_ZhCh} 
\end{equation}
with  $\sigma, m_i$, $i=1,2,3,4$, arbitrary parameters.
The first term 
%\begin{equation}
$ \sigma \left( J_3^2 - \beta^2 J_1^2 \right)$
% \label{eq:W_Euler} 
%\end{equation}
is equivalent (up to an affine transformation) to the Hamiltonian of the quantum Euler top (see \cite{PW_Heun, T}).
The complete operator \eqref{eq:W_ZhCh} corresponds to the Hamiltonian of a Euler top
 %(quadratic terms in \eqref{eq:W_ZhCh}) 
 with additional \textquotedblleft magnetic" interactions accounted for  by the linear terms in \eqref{eq:W_ZhCh}. In classical mechanics this is the Zhukovsky-Volterra gyrostat \cite{Basak}, \cite{LOZ}. Note moreover that a similar Hamiltonian was exploited to describe spin systems with anisotropy \cite{ZU}.
We thus see that the Heun operator pencil on the $\su$ algebra is equivalent to the Hamiltonian of a generalized Euler top (or quantum Zhukovskii-Volterra gyrostat). 

\section{The case of the Lie algebra $\sun$}

\subsection{Specializations\label{sec:spec3}}

The situation where $\sun$ is the underlying algebra will have, not surprisingly, close similarities with the picture in the $\su$ case; the main differences will come 
from the richer orbit structure. Let us indeed explain how the choices for the bispectral pair $(X,Y)$ can be standardized. Throughout this section, the elements $J_\kappa$
with $\kappa = \:1,2,3$ obey the relations corresponding to the lower signs in \eqref{eq:crsu} and \eqref{eq:csu}.
So as to relate to the unitary representations
of $\sun$ and to have one discrete spectrum in play, we shall again take $Y$ to be the compact generator $J_3$. Following a reasoning similar to the one 
given in the last section, we observe that a generic element  $ l_1J_1+l_2J_2+l_3J_3$ can be transformed into $\alpha J_3 +\beta J_1$ under a rotation while not changing $Y$.
The difference here is that the adjoint action of $\sun$ preserves the non-definite form $l^2=l_1^2+l_2^2-l_3^2$. This will hence require that $\beta^2-\alpha^2 =l^2$. 
There will therefore be three distinct cases according to whether $l^2$ is negative, positive or zero; elements in these classes are respectively said to be elliptic,
hyperbolic and parabolic.

The upshot of this is that the bispectral operators $X$ and $Y$ in the $\sun$ case will have exactly the same form as those of the $\su$ case, namely
\begin{equation}
 X=\alpha J_3 +\beta J_1 \quad, \quad  Y= J_3 \label{eq:XYsun}
\end{equation}
but will fall into the three categories where allowing for scaling, $\alpha$ and $\beta$ will be parametrized as follows:
%\begin{itemize}
%\item  \qquad elliptic \: \;  \ \qquad $l^2 <0$ \qquad $\alpha=ch\theta$ \; \qquad $\beta=-sh\theta$
%\item \qquad hyperbolic \qquad $l^2>0$ \qquad $\alpha =-cos\phi$ \qquad $\beta = 1$
%\item \qquad parabolic \;  \qquad $l^2=0$ \qquad $\alpha=1$ \qquad \ \qquad $\beta=1.$
%\end{itemize}
%$0<\phi<\p$
%\begin{eqnarray}
%&&\text{elliptic} \:\;\;\; \;\; \;\qquad l^2 <0 \: \qquad \alpha=ch\theta \qquad \beta=-sh\theta \label{eq:ell}\\
%&&\text{hyperbolic}\:\: \qquad l^2 >0  \qquad \alpha =-cos\phi \qquad \beta = 1 \label{eq:hyp} \\
%&& \text{parabolic}\: \; \; \: \qquad l^2=0 \qquad \alpha=1 \qquad \ \qquad \beta=1 \label{eq:para}\
%end{eqnarray}
\begin{align}
&\text{elliptic}&&l^2 <0 && \alpha=\cosh\theta && \beta=-\sinh\theta \label{eq:ell}\\
&\text{hyperbolic}&& l^2 >0  && \alpha =-\cos\phi && \beta = 1 \label{eq:hyp} \\
& \text{parabolic}&& l^2=0 && \alpha=1 && \beta=1 \label{eq:para}\
\end{align}
with  $ 0< \phi <\pi$.

Owing to the fact that $X$ and $Y$ have the same general expressions in both the $\su$ and $\sun$ cases many formulas will be almost identical to those for $\su$ as
we determine the Heun algebras associated to $\sun$. 
$X$ and $Y$ will again obey the relations \eqref{eq:ZX} and \eqref{eq:YZ} of the algebra $\cA$ with the only change with respect to \eqref{eq:par} being
\begin{equation}
c_2=\beta^2-\alpha^2. \label{eq:coefc_1}
\end{equation}
The Casimir operator $C$ of $\cA$ will be related exactly as in \eqref{eq:Cas} to the Casimir element $\mathfrak{c}$ of $\sun$ given by the lower part of \eqref{eq:csu}.
We may point out that the relations \eqref{eq:conj} remain true for the elliptic case. 
Always keeping the definition \eqref{eq:W} for $W$, Proposition \ref{propHeun} readily translates to the $\sun$ case.
\begin{prop}\label{pr:heunc}
The elements $X$ and $W$ associated to $\sun$ as per \eqref{eq:XYsun} satisfy the relations of the Heun algebra of Lie type with the coefficients $x_s$ and $y_s$ ,
$ s=0, 1, 2, 3$ as in \eqref{eq:Ha1} -  \eqref{eq:Ha2} and with $\alpha$ given by \eqref{eq:ell}, \eqref{eq:hyp}, \eqref{eq:para} for each of the possible classes and $c_2$
given by \eqref{eq:coefc_1}.
\end{prop}
The central element $\Omega$ of the corresponding Heun algebras is again given by \eqref{eq:omega} with the appropriate $\alpha$ and $c_2$.

\subsection{Connections to the Hahn algebra \label{sec:hahn}}

Having found for $\sun$, three Heun algebras corresponding to whether $Y$ is of elliptic, hyperbolic or parabolic type, we now observe that these
algebras are isomorphic to the Hahn algebra. 
Indeed, the relation \eqref{eq:mn} for the parameters $\mu$ and $\nu$ has real solutions that read in light of \eqref{eq:coefc_1}:
\begin{equation}
 \nu=-\frac{2r_2}{\alpha\pm \beta} \quad\text{and} \qquad \mu=-\frac{r_4}{\alpha\pm \beta} - r_3.
\end{equation}
For the parabolic case $\alpha=\beta=1$, only the solution with the upper sign is permitted. The operators $X$ and $\overline{W}$
therefore satisfy the Hahn algebra with the parameters 
\begin{eqnarray}
 && \overline x_0 =(\alpha^2-\beta^2)r_5 \ , \qquad \overline x_1=\pm\beta r_4 \ , \qquad \overline x_2=\pm 2\beta r_2 \ , \qquad \overline x_3=\beta^2-\alpha^2 \ , \label{eq:cc1}\\
 &&\overline y_0= \mp \beta r_4 r_5-2r_2r_4C , \qquad \overline y_1=\beta^2(r_2^2-r_1^2)\mp4\beta r_2r_5 -4r_2^2 C\ . \label{eq:cc2}
\end{eqnarray}

\section{The case of the harmonic oscillator algebra \label{sec:osc}}

We shall focus in this section on the situation where the algebra $\cA$ of Definition \ref{defn1} is isomorphic to the harmonic oscillator algebra $\mathfrak{ho}$ \eqref{eq:reho}.
Consider for the bispectral pair $(X,Y)$:
\begin{equation}
X=N+\chi(A+A^{\dag}) +\chi^2\; I \quad , \quad Y=N \ ,\label{eq:ho} 
\end{equation}
where $\chi$ is a real constant. Here again, it is straightforward to check that the relations \eqref{eq:ZX} and \eqref{eq:YZ} of the algebra $\cA$ are satisfied with
\begin{equation}
b=1 \qquad c_1=c_2=-1 \quad \quad d_1=d_2=\chi^2
\end{equation}
and one readily observes that $X, Y, Z$  and  $I$ also form a basis for $\mathfrak{ho}$. The Casimir $C$ of $\cA$ becomes:
\begin{equation}
 C=2\chi^2(X+Y)+ \{X,Y\}-X^2-Y^2+Z^2= \chi^2 [4 \mathfrak{k}+\chi^2-2]
\end{equation}
with $\mathfrak{k}=N-A^{\dag}A$.
We now bring on the associated algebraic Heun operator $W$ given in \eqref{eq:W} as the generic bilinear expression in these $X$ and $Y$.
Anew, it is seen that the generator $Y$ is obtained from $X$ by an automorphism:
\begin{equation}
Y=U(\chi)XU(\chi)^{-1} \qquad \text{with} \qquad U(\chi)=e^{\chi(A-A^{\dag})}.
\end{equation}
Recalling from the argument given in Section 2 that the Heun algebras generated by the pairs $(X,W)$ and $(Y,W)$ are 
isomorphic, we shall only concentrate on the former with the following result directly obtained from Proposition \ref{prop1}.
\begin{prop}\label{pr:charlier} The elements $X$ and $W$ associated to $\mathfrak{ho}$ as per \eqref{eq:ho} and \eqref{eq:W} satisfy the relations of the Heun algebra of Lie type with
\begin{eqnarray}
&&x_0= \chi^2 r_4+r_5 \ ,   \qquad x_1=2\chi^2 r_2+r_3+r_4\ , \qquad x_2=2r_2\ , \qquad x_3=-1\ ,\\
&&y_0=\chi^2(2r_1^2-2r_2^2-r_4^2- r_3r_4+2r_2r_5)-(r_3+r_4)r_5-2r_2r_4 C\ ,\\
&&y_1=4\chi^2 r_2  (2r_4-r_3)-(r_3+r_4)^2-4r_2r_5 -4r_2^2C\ ,\\
&&y_2=12\chi^2 r_2^2-6r_2(r_3+r_4)\ , \qquad y_3=-8r_2^2\ .
\end{eqnarray}
\end{prop}

With these parameters, the central element $\Omega$ of $\cH$ is 
\begin{equation}
 \Omega=(r_4^2+r_2^2-r_1^2)C+\chi^4(r_1^2-r_2^2)-2\chi^2r_4(r_2-r_5)+2r_2r_5+r_5^2 \;.
\end{equation}
In this case, looking at the relations \eqref{eq:mn} for $\mu$ and $\nu$ one finds that there are no solutions.
The Heun algebra associated to the oscillator algebra is therefore not isomorphic to the Hahn algebra.

\section{Realizations in terms of difference operators}

It is well known that various families of orthogonal polynomials are related to the Lie algebras that we have considered so far \cite{GZ2,FLV}.
We shall review these results in this section by considering different models for $\su$, $\sun$ and $\mathfrak{ho}$ and by
showing that the associated operators $Y$ are realized within the corresponding representations as the operators of which
the polynomials are eigenfunctions. Since $X$ will always be multiplication by the variable, this will allow us to say that the
Heun algebras of Lie type that have been identified are the Heun algebras of Krawtchouk, Meixner, Meixner-Pollaczek,
Laguerre and Charlier type.

\subsection{$\su$}

There is a model of $\su$ in terms of the difference operators $T^\pm$
\begin{equation}
T^\pm f(x) = f(x\pm 1).\label{eq:shr}
\end{equation}
It has the generators given as follows:
\begin{eqnarray}
 J_3&=&\sin^2(\theta/2)\ (x-N)T^+  +\cos(\theta)\ (x-N/2)  -\cos^2(\theta/2)\ xT^-\ , \label{eq:J_3}\\
 J_-&=& \frac{\sin\theta}{2}\left(  (x-N)T^+ -2x+N+xT^-   \right) \ , \\
 J_+&=& \frac{\sin\theta}{2}\left( \tan^2(\theta/2) (N-x)T^+ -2x+N-\cot^2(\theta/2) xT^-   \right) \ ,
\end{eqnarray}
with $J_1=\frac{1}{2}(J_+ + J_-)$ and $J_2=-\frac{i}{2}(J_+ - J_-)$.
As per Section \ref{sec:spec2}, we have the bispectral operators $X$ and $Y$ given by
\begin{eqnarray}
 X=\cos(\theta)J_3-\sin(\theta)J_1 =x-\frac{N}{2} \quad\text{and}\quad Y=J_3\ .
\end{eqnarray}
One readily recognizes in view of \eqref{eq:J_3}, that $Y$ becomes the difference operator of the Krawtchouk polynomials $K_n(x;\sin^2(\theta/2),N)$ \cite{KLS}:
\begin{equation}
 Y K_n(x;\sin^2(\theta/2),N)=\left(n-\frac{N}{2}\right)K_n(x;\sin^2(\theta/2),N)\ .
\end{equation}
The operators $J_\pm=J_1\pm iJ_2$ are moreover lowering and raising operators for these polynomials:
\begin{eqnarray}
 &&J_-K_n(x;\sin^2(\theta/2),N)= n\cot(\theta/2) K_{n-1}(x;\sin^2(\theta/2),N)\ , \\
 &&J_+K_n(x;\sin^2(\theta/2),N)= (N-n)\tan(\theta/2) K_{n+1}(x;\sin^2(\theta/2),N) \ .
\end{eqnarray}
The associated operator $W$
\begin{eqnarray}
W&=& \sin^2(\frac{\theta}{2}) (N-x)(Nr_2+r_1-r_2-r_4-2r_2x)T^+ +2r_2\cos(\theta)x^2+\rho_4 x +\rho_5\nonumber\\
&&+\cos^2(\frac{\theta}{2})x(Nr_2-r_1+r_2-r_4-2r_2x)T^-
\end{eqnarray}
is hence the Heun-Krawtchouk operator. 
The parameters $\rho_4$ and $\rho_5$ are given by
\begin{eqnarray}
 \rho_4=\cos(\theta)(r_4-2Nr_2)+r_3\ , \qquad
\rho_5=\frac{1}{2}\cos(\theta)N(Nr_2-r_4)-\frac{Nr_3}{2}+r_5\ .
\end{eqnarray}
The algebra given in Proposition \ref{propHeun} can thus appropriately be called the Heun-Krawtchouk algebra since its generators can be 
realized using the bispectral operators of the Krawtchouk polynomials.

\subsection{$\sun$} Three models of the positive discrete series representation will be presented in correspondance with the situations
where $Y$ is of elliptic, hyperbolic and parabolic type. The orthogonal polynomials of Meixner, Meixner-Pollaczek and Laguerre will
arise correspondingly \cite{MR, KV}. Furthermore, owing to the fact that the associated Heun algebras are isomorphic to the Hahn one,
it will be recorded that the respective Heun operators $W$ can be transformed into the Hahn, Continuous Hahn and Jacobi operators.

\subsubsection{Elliptic case}

The generators of $\sun$ can be realized as follows:
\begin{eqnarray}
 J_3&=&- (x+\kappa) \sinh^2(\theta/2) T^+ +(x+\kappa/2)\cosh(\theta)    -\cosh^2(\theta/2) xT^- \\
 J_-&=& \frac{-\sinh\theta}{2}\left(  (x+\kappa)T^+ -2x-\kappa + xT^-   \right) \\
 J_+&=& \frac{-\sinh\theta}{2}\left( \tanh^2(\theta/2) (x+\kappa)T^+  -2x-\kappa + \coth^2(\theta/2) xT^-   \right) \ .
\end{eqnarray}
in terms of the shift operators \eqref{eq:shr}.
As per Section \ref{sec:spec3}, we introduce the following bispectral operators
\begin{eqnarray}
 X=\cosh(\theta)J_3-\sinh(\theta)J_1 =x+\frac{\kappa}{2} \quad\text{and}\quad Y=J_3\ .
\end{eqnarray}
$Y$ is then identified as the difference operator of the Meixner polynomials $M_n(x;\kappa,\tanh^2(\theta/2))$ \cite{KLS}:
\begin{equation}
 Y  M_n(x;\kappa,\tanh^2(\theta/2))=\left(\frac{\kappa}{2}+n\right) M_n(x;\kappa,\tanh^2(\theta/2)) \ 
\end{equation}
and the operators $J_+$ and $J_-$ are ladder operators for these polynomials:
\begin{eqnarray}
 &&J_-M_n(x;\kappa,\tanh^2(\theta/2))= n\coth(\theta/2) M_{n-1}(x;\kappa,\tanh^2(\theta/2)) \\
 &&J_+M_n(x;\kappa,\tanh^2(\theta/2))= (\kappa+n)\tanh(\theta/2) M_{n+1}(x;\kappa,\tanh^2(\theta/2)) \ .
\end{eqnarray}
The associated operator $W$
\begin{eqnarray}
W&=&-(x+\kappa)(2r_2 x+r_2\kappa-r_1+r_2+r_4)\sinh^2(\theta/2)T^+  +2r_2\cosh(\theta)x^2+\rho_4 x+\rho_5\nonumber\\
&&+ x(r_2 x+r_2 \kappa+r_1-r_2+r_4)\cosh^2(\theta/2)T^- 
\end{eqnarray}
is hence the Heun-Meixner operator. We have introduced
\begin{eqnarray}
 \rho_4=\cosh(\theta)(2r_2\kappa+r_4)+r_3\ , \qquad \rho_5=  \frac{\cosh(\theta)\kappa}{2}(r_2\kappa+r_4) +\frac{r_3\kappa}{2}+r_5\ .
\end{eqnarray}
The algebra given in Proposition \ref{pr:heunc} may thus be called the Heun-Meixner algebra since it
can be realized using the bispectral operators of the Meixner polynomials.
We know that this Heun-Meixner algebra is isomorphic to the Hahn one. This can be underscored within
the present $\sun$ representation by recovering the Hahn operator from W
%In this case, we can also obtain a model of the Hahn algebra with the  operators $X$ and $\overline W$ as generators
%by 
recalling the observations made in Section \ref{sec:hahn}. 
%To simplify the results, we 
Consider to simplify, the conjugated operator $\widetilde W=\frac{1}{\tanh^x(\theta/2)}\overline W\tanh^x(\theta/2)$ which reads as follows
\begin{eqnarray}
 \widetilde W&=&\frac{1}{\tanh^x(\theta/2)} W\tanh^x(\theta/2) -(r_4 e^{-\theta} + r_3)X-2r_2e^{-\theta} X^2\nonumber\\
 &=& -r_2\sinh(\theta)\Big[x(x+\frac{r_2\kappa+r_1-r_2+r_4}{2r_2})T^-+(\kappa+x)(x+\frac{r_2\kappa-r_1+r_2+r_4}{r_2})T^+\nonumber\\
 &&\hspace{2cm}-(2x+\kappa)(x+\frac{r_2\kappa+r_4}{2r_2})\Big]+r_5 \label{eq:WE}
\end{eqnarray}
We readily recognize 
the difference operator of the Hahn polynomials in the square bracket of the second line of \eqref{eq:WE}.

\subsubsection{Hyperbolic case}

The operators
\begin{eqnarray}
 J_3&=&\frac{1}{2i\sin\phi}\left(e^{i\phi}(\lambda-ix)\mathcal{T}^+ +2ix\cos(\phi) -e^{-i\phi}(\lambda+ix)\mathcal{T}^-\right)\\
 J_\pm&=&\frac{e^{i\phi}}{2i\sin\phi}\left(e^{\pm i\phi}(\lambda-ix)\mathcal{T}^+ +2ix -e^{\mp i\phi}(\lambda+ix)\mathcal{T}^-\right)
\end{eqnarray}
satisfy the $\sun$ commutation relations with the shift operators $\mathcal{T}^\pm$ defined by $\mathcal{T}^\pm f(x)=f(x\pm i)$.
With $X$ and $Y$ given according to Section \ref{sec:spec3} by
\begin{eqnarray}
X&=& - \cos(\phi)J_3 +J_1=x\sin(\phi)\ , \qquad Y=J_3,
\end{eqnarray}
we see that $Y$ is realized in this instance as the difference operator of the Meixner-Pollaczek polynomials $ P_n^{(\lambda)}(x;\phi)$ \cite{KLS}:
\begin{equation}
 Y P_n^{(\lambda)}(x;\phi)=(n+\lambda) P_n^{(\lambda)}(x;\phi)\ ,
\end{equation}
while $J_+$ and $J_-$ act as raising and lowering operators:
\begin{eqnarray}
\hspace{1cm} J_+  P_n^{(\lambda)}(x;\phi)=(n+1)P_{n+1}^{(\lambda)}(x;\phi)\ , \qquad
 J_-  P_n^{(\lambda)}(x;\phi)=(2\lambda+n-1)P_{n-1}^{(\lambda)}(x;\phi)\ .
\end{eqnarray}
The associated operator $W$
\begin{eqnarray}
W&=&(ix-\lambda)\left(ir_2x+ \frac{r_1-r_2}{2}+\frac{ir_4}{2\sin(\phi)} \right)\mathcal{T}^+ +2r_2\cos(\phi)x^2 +\rho_4 x +r_5 \nonumber\\
&&+e^{-i\phi}(\lambda+ix)\left( ir_2x - \frac{r_1-r_2}{2} +\frac{ir_4}{2\sin(\phi)} \right)\mathcal{T}^-
\end{eqnarray}
is hence the Heun-Meixner-Pollaczek operator. 
Here $\rho_4=\sin(\phi)r_3 +\cot(\phi) r_4$.
In this case the algebra of Proposition \ref{pr:heunc} is really the Heun-Meixner-Pollaczek algebra.
In this realization, upon scaling and conjugating $\overline{W}$ according to $\widetilde W=\frac{1}{r_2} e^{i\phi x} \overline{W} e^{-i\phi x}$ one finds
\begin{eqnarray}
 \widetilde{W}&=&(\lambda-ix)\left( \frac{r_2-r_1}{2r_2}-\frac{ir_4}{2r_2\sin(\phi)} -ix\right)\mathcal{T}^+ -2x^2-\frac{x  r_4}{r_2\sin(\phi)}+\frac{r_5}{r_2}\nonumber\\
 &&+(\lambda+ix)\left( \frac{r_2-r_1}{2r_2} +\frac{ir_4}{2r_2\sin(\phi) } +ix \right)\mathcal{T}^-.
\end{eqnarray}
This operator is the difference operator that is diagonalized by the continuous Hahn polynomials 
$(-1)^{ix}p_n(x;\lambda, \frac{r_2-r_1}{2r_2} +\frac{ir_4}{2r_2\sin(\phi) },\lambda ,\frac{r_2-r_1}{2r_2}-\frac{ir_4}{2r_2\sin(\phi)})$ with 
 eigenvalues $-n^2-(2n+1)\lambda+ \frac{(\lambda+n)r_1+r_5}{r_2}$.

\subsubsection{Parabolic case}

A model of $\sun$ in terms of differential operator is given by
\begin{eqnarray}
 J_3&=&-x\frac{d^2}{dx^2}-(1+a-x)\frac{d}{dx}+\frac{1+a}{2}\\
 J_-&=&   x \frac{d^2}{dx^2} +(1+a)\frac{d}{dx} \\
 J_+&=&  x \frac{d^2}{dx^2} +(1+a-2x)\frac{d}{dx}+x-1-a
\end{eqnarray}
where $a$ is a free parameter associated to this realization. For the parabolic case, the bispectal operators $X$ and $Y$ were taken to be 
\begin{eqnarray}
 X=J_3+J_1=\frac{1}{2}\ x \quad\text{and}\quad Y=J_3\ 
\end{eqnarray}
in Section \ref{sec:spec3}.
It follows that $Y$ can here be identified with the difference operator of the Laguerre polynomials $ L_n^{(a)}(x)$ \cite{KLS}:
\begin{equation}
 Y L_n^{(a)}(x)=(n+(a+1)/2) L_n^{(a)}(x)\ .
\end{equation}
One also gets 
\begin{eqnarray}
&& J_+L_n^{(a)}(x)=-(n+1)L_{n+1}^{(a)}(x)\ , \qquad  J_-L_n^{(a)}(x)=-(a+n)L_{n-1}^{(a)}(x)\ .
\end{eqnarray}
The associated operator $W$
\begin{equation}
W=-x(r_2x+r_4) \frac{d^2}{dx^2}+(r_2 x^2+(r_1-2r_2-r_2a+r_4)x-r_4(1+a))\frac{d}{dx}
+ \rho_4 x +\rho_5  
\end{equation}
is hence the Heun-Laguerre operator. 
We have defined
\begin{eqnarray}
 \rho_4=\frac{1}{2}( r_3 - r_1 + (a+2)r_2 )\quad\text{and}\qquad
 \rho_5=r_5 + \frac{1}{2}(r_1-r_2+r_4  )(a+1) \ .
\end{eqnarray}
In this case, the algebra given in Proposition \ref{pr:heunc} should hence be referred to as  the Heun-Laguerre algebra since 
we have a realization of it based on the bispectral operators of the Laguerre polynomials.
It can be also showed that the eigenvalue equation $Wf(x)=\lambda f(x)$ becomes the confluent Heun equation \cite{K}, \cite{Ro}. Let us remark that
this connection between the confluent Heun equation and the Heun operators associated to the Laguerre differential equation had already been pointed out in \cite{GVZ2}.
Finally note that the conjugated operator $\widetilde {W}=e^{-x/2}\overline We^{x/2}$ 
\begin{eqnarray}
 \widetilde W&=&e^{-x/2}We^{x/2}-(r_4/2 + r_3)X-r_2 X^2\nonumber\\
 &=& -x(r_2x+r_4) \frac{d^2}{dx^2}+( (r_1-2r_2a-r_2)  x-r_4(1+a))\frac{d}{dx}
+\frac{1}{2}(r_1-r_2)(a+1)+r_5  
\end{eqnarray}
is recognized to be the Jacobi differential operator 
which together with $X$ generates the Hahn algebra. Indeed, $\widetilde W$ is diagonalized by the Jacobi polynomial
$P_n^{(a,-r_1/r_2)}\left(\frac{2r_2x}{r_4}+1\right)$ with eigenvalues $ (r_1-r_2)(n+\frac{1}{2}(a+1) )-r_2 n(n+a)+ r_5$.

\subsection{Harmonic oscillator algebra $\mathfrak{ho}$}

The oscillator algebra can also be realized as follows in terms of the shift operators \eqref{eq:shr}: 
\begin{equation}
 N=-xT^-+x+\chi^2-\chi^2 T^+  \ , \quad A^\dagger= -\chi +\frac{1}{\chi} xT^- \ , \quad A=\chi(T^+-1) \ .
\end{equation}
As in Section \ref{sec:osc}, the bispectral pair $(X,Y)$ is
\begin{eqnarray}
 X=N+\chi(A+A^\dagger)+\chi^2=x\ , \qquad Y=N \ .
\end{eqnarray}
The operator $Y$ is the difference operator of the Charlier polynomial $C_n(x,\chi^2)$ \cite{KLS}:
\begin{equation}
 YC_n(x,\chi^2)=nC_n(x,\chi^2)\ ,
\end{equation}
and $A^\dagger$ and $A$ are their raising and lowering operators:
\begin{equation}
 A^\dagger C_n(x,\chi^2)=-\chi C_{n+1}(x,\chi^2)\ , \qquad A C_n(x,\chi^2)=-\frac{n}{\chi} C_{n-1}(x,\chi^2)\ .
\end{equation}
The associated operator $W$
\begin{equation}
W=x(r_2-r_4-r_1-2r_2 x)T^- +\chi^2( r_1   -r_2-r_4  -2 r_2 x)T^+ +2r_2x^2  +(2r_2\chi^2+r_3+r_4)x+r_4\chi^2+r_5
\end{equation}
is hence the Heun-Charlier operator.
The algebra given in Proposition \ref{pr:charlier} should therefore be called the Heun-Charlier algebra since 
it is realized with $X=x$ and $W$ constructed from $x$ and the Charlier operator.

\section{Concluding remarks}
The results presented here provide a comprehensive picture of the Heun operators and algebras associated to orthogonal
polynomials that admit a Lie theoretical interpretation. This complements the previous studies of the Jacobi and Hahn polynomials
from this Heun angle which led respectively to the description of the standard Heun differential operator and its discrete version
within the framework of bispectral problems. Missing is the parallel treatment of the Racah polynomials especially since the $q \rightarrow 1$ limit
of the Askey-Wilson case is known to be delicate. It would also be quite pertinent to carry on with the q-Askey tableau and to work out the quantum algebraic analog of the study
offered in the present article. It is known that Heun operators have applications in time and band limiting problems, it would be of interest to look
at applications in this direction in particular. We hope to report on these questions in the near future.

\
\medskip

 \textbf{Acknowledgments:} N.C. is gratefully holding a CRM--Simons professorship. 
 The research of L.V. is supported in part by a Natural Science and Engineering Council (NSERC) of Canada
discovery grant and that of A.Z. by the National Science Foundation of China (Grant No. 11711015).

\end{document}